\documentclass{article}

\usepackage[final,nonatbib]{neurips_2020_tda}

\usepackage[utf8]{inputenc} 
\usepackage[T1]{fontenc}    
\usepackage{amsfonts}       
\usepackage{booktabs}       
\usepackage{hyperref}       
\usepackage{url}            
\usepackage{microtype}      
\usepackage{nicefrac}       
\usepackage{paralist}       
\usepackage{siunitx}        
\usepackage{amsmath,amssymb}
\usepackage{mymacros}
\usepackage{graphicx}

\title{Multidimensional Persistence Module Classification via Lattice-Theoretic Convolutions}

\author{%
  Hans Riess \\
  University of Pennsylvania\\
  \texttt{hmr@seas.upenn.edu} \\
  \And
  Jakob Hansen \\
  BlueLight AI \\
  \texttt{jakob@jakobhansen.org}\\
  \And
  Robert Ghrist \\
  University of Pennsylvnaia \\
  \texttt{ghrist@seas.upenn.edu}
}

\begin{document}

\maketitle

\begin{abstract}
 Multiparameter persistent homology has been largely neglected as an input to
 machine learning algorithms.
 We consider the use of lattice-based convolutional neural network layers as a 
 tool for the analysis of features arising from multiparameter persistence
 modules. We find that these show promise as an alternative to convolutions for
 the classification of multidimensional persistence modules.
\end{abstract}

\section{Introduction}

Persistent homology has the ability to discern both the global
topology~\cite{ghrist_barcodes:_2008} and local
geometry~\cite{bubenik_persistent_2020} of finite metric spaces (e.g. embedded
weighted graphs, point clouds in $\R^d$) making it a befitting feature for the
purposes of training a neural network. Single-dimensional homological
persistence has drawn recent attention in deep learning~\cite{hofer_deep_2017,
pun_persistent-homology-based_2018, bruel-gabrielsson_topology_2020}. This is,
in part, due to a wide range of efficient software
libraries~\cite{otter_roadmap_2017, henselman_matroid_2017, bauer_ripser:_2019}
for computing persistent homology, 
as well as a growing cookbook of recipes 
for featurizing barcodes from single dimensional persistence, including
persistence images~\cite{adams_persistence_2017}, persistence
landscapes~\cite{bubenik_statistical_2015}, and more exotic methods
\cite{kalisnik_tropical_2019}.

Multidimensional persistence generalizes single-dimensional persistent homology
in order to tackle filtrations parameterized in multiple dimensions.
Unfortunately, there is no complete
compact barcode-like characterization of multidimensional persistence
modules~\cite{carlsson_theory_2009}. We must make do with incomplete invariants.
There are various algebraic invariants; in this paper we will use the
\textit{Hilbert function} and the \textit{multi-graded Betti numbers}\footnote{\textit{Caveat lector:} the multi-graded (algebraic) Betti numbers need not be confused with the topological Betti numbers, i.e.~the rank of homology.}, both of
which are $\mathbb{N}$-valued functions on the parameter space.
The Hilbert function is nothing more than the pointwise
(topological) Betti numbers, and the multi-graded Betti
numbers have a geometric interpretation in terms of births and deaths~\cite{knudson_refinement_2008}.
 The \textit{rank
 invariant}, another invariant, is shown to be complete in the case of single
 dimensional persistence \cite{carlsson_theory_2009}; due to its difficulty to
 compute it is not considered here.

Multidimensional persistence suffers two hindrances to its usefulness in
machine learning. First, software for computing multidimensional persistence is
scarce. To the authors' knowledge, RIVET is the only available software for
computing multidimensional persistence~\cite{lesnick_interactive_2015}; RIVET
specializes to $2$-dimensional persistence and is focused on interactive
visualization rather than machine-interpretable output.
Second, 
there has been only very recent and preliminary activity \cite{vipond2020multiparameter,carriere2020multiparameter} focusing on extraction
of suitable features for multidimensional persistent homology. A full-fledged deep learning pipeline
for multidimensional persistence waits in the wings.

We hope to ignite further interest in filling both of these gaps. In this paper, we
propose a naive featurization of multidimensional persistence modules based on
the aforementioned invariants, and design an architecture for classifying these
persistence modules. This architecture employs a lattice-theoretic notion of
convolution, thereby respecting the order relation of the parameters of the
persistence module. We implement our model and compare the performance of our
proposed lattice-convolutional architecture with a (simplified) standard
convolutional architecture.

\paragraph{Related work}
Very recent advances have been made to featurize multidimensional persistence modules via
landscapes \cite{vipond2020multiparameter} and images \cite{carriere2020multiparameter}.
Our naive approach at featurizing $2$-dimensional persistence modules has the advantage
of being readily computable with existing software \cite{lesnick_interactive_2015}.
Persistence modules supported on lattices are a popular object of study as of late.
By leveraging M\"{o}bius inversion, these persistence modules are shown to generate a stable persistence diagram \cite{mccleary2020edit} as well as factor through a convenient category with pseudo-inverses \cite{krishnan2020invertibility}.

\section{Backgound}

Space constraints require that this section be laconic. For a primer on
persistent homology, see~\cite{ghrist_barcodes:_2008, carlsson_topology_2009};
for multiparameter persistent homology, see~\cite{lesnick_interactive_2015}. An
introduction to lattices may be found in~\cite{davey2002introduction}.

\subsection{Rips complexes and persistent homology}
Let $(\mathcal M,d)$ be a finite metric space. The \introduce{Vietoris-Rips
  complex} of $\mathcal M$ at scale $r$ is the abstract simplicial complex
$\text{Rips}_r(\mathcal M)$ whose simplices are subsets of $\mathcal M$ of
diameter at most $r$. There is a natural inclusion $\text{Rips}_r(\mathcal M)
\to \text{Rips}_{r'}(\mathcal M)$ for $r \leq r'$.

Applying the simplicial homology functor (with coefficients in a field $k$)
$H_i$ to $\text{Rips}_r(\mathcal M)$ produces a sequence of vector spaces
$PH_i(r)$. The inclusions $\text{Rips}_r(\mathcal M) \to
\text{Rips}_{r'}(\mathcal M)$ induce maps $PH_i(r) \to PH_i(r')$, producing the
data of a \introduce{persistence module}. This structure can be compactly
described as a functor from $\R$, viewed as a category via its standard order
structure, to the category $\Vect_k$ of vector spaces over $k$. The simplicity
of the category $\R$ gives these persistence modules simple structure: they
decompose as direct sums of interval modules $I_{[a,b)}$, which have
$I_{[a,b)}(r) = k$ for $a \leq r < b$ and zero otherwise. The maps are the
identity where possible and the zero map otherwise.


\subsection{Multiparameter persistence}
The Rips construction produces a filtration of simplicial complexes from a
finite metric space; it is natural to consider the behavior of the homology
functor over a pair of coherent filtrations. Consider
a finite metric space $(\mathcal M, d)$ and a filtration
function $\rho:\mathcal M \to \R$. This data specifies a bifiltration of
simplicial complexes given by
\[\mathbb{X}_{r,t} = \text{Rips}_r\{x \in \mathcal M \mid \rho(x) \leq t\}.\]
There is a natural inclusion $\mathbb{X}_{r,t} \hookrightarrow
\mathbb{X}_{r',t'}$ whenever $(r,t) \leq (r',t')$ in the lattice $\R \times \R$.
Composing with the homology functor produces a 2-parameter persistence module
\[PH_{i}: \R_+ \times \R \to \Vect; \quad (r,t) \mapsto H_i(X_{r,t}).\]


While there does not exist a complete discrete invariant for  $PH_i$, we can extract
meaningful features. Two particularly informative types of features are 
 the \introduce{Hilbert function}
  \[\text{Hilb}: \R_+ \times \R \to \Z_+;\quad (r,t) \mapsto \dim(PH_i(r,t)),\] and
the \introduce{multi-graded Betti numbers} $\xi_j: \R_+ \times \R \to \Z_+$, for $j =
0,1,2$. For $PH_i$ as above, the Hilbert function counts the number of
connected components ($i=0$), cycles ($i=1$), or higher dimensional voids ($i >
1$) of the complex $\mathbb{X}_{r,t}$ at each $(r,t) \in \R_+\times \R$. The multi-graded
Betti numbers, on the other hand, capture information about births
and deaths of homology classes.\footnote{ $\xi_j(r,t)$ counts the rank of the $j$th term of the free resolution of $PH_i$ at grading $(r,t)$ \cite{lesnick_interactive_2015}.}

\subsection{Lattice-theoretic signal processing}\label{sec:latticeconv}

Classical signal processing proceeds by constructing filters for space- or time-indexed signals using convolutional operators. These implicitly rely on the
algebraic properties of the underlying space, in particular the existence of
well-behaved translation operators. This fact is exploited in the general
framework of algebraic signal processing~\cite{puschel_algebraic_2008} and
extended to more general domains by the theory of graph signal
processing~\cite{ortega_graph_2018}. We here describe a similar extension to
signals on a finite lattice proposed in~\cite{pueschel_discrete_2019}.

A \introduce{lattice} $L$ is a partially ordered set in which every pair of
elements $x,y$ has a greatest lower bound (the \introduce{meet} $x \wedge y$)
and a least upper bound (the \introduce{join} $x \vee y$). These operations
and their properties produce an algebraic characterization of lattices; the
ordering can be recovered from the algebra and vice versa. The key insight of~\cite{pueschel_discrete_2019} is that the meet and
join operations on a lattice define two ``shift operators'' that can be
exploited to define convolutional filters for signals on a lattice.

That is, for two signals $f,g: L \to \R$, where $L$ is a lattice, we define
\[(f *_\wedge g)(x) = \sum_{a \in L}f(x\wedge a)g(a) \;\;\text{and}\;\; (f *_\vee g)(x) = \sum_{a \in L}f(x\vee a)g(a).\]

A nice class of examples of lattices are given by the sets $\R^n$
viewed as partially ordered sets, with the ordering $(x_1,\ldots,x_n) \leq
(y_1,\ldots,y_n)$ whenever $x_i \leq y_i$ for all $1 \leq i \leq n$. These are,
of course, the indexing sets for persistent homology, suggesting that lattice
convolutions over $\R^n$ or its finite sublattices may be useful in processing
data coming from multiparameter persistence computations. 

\section{Lattice Convolutional Neural Networks}\label{sec:latticeCNN}

Convolutions over $\R^2$ (with its abelian group structure)
have served as an easily parameterized and efficient set of linear operations
adapted to the structure of images. Their extreme utility in computer vision
problems is owed to the translation equivariance properties of images: humans
naturally recognize an image translated via an additive reparameterization as
equivalent to the original.

The data of a multidimensional persistence module is also indexed by $\R^n$ or a
regular finite subset thereof, but its natural algebraic structure is not that
of an abelian group. Rather, with its partial order structure, the indexing set
is a lattice. In processing signals associated with a persistence module, it
may be useful to take this structure into account rather than imposing the
abelian group structure implied by standard convolutions.

To this end, we construct a lattice convolution-based neural network layer
suitable for use with features originating from multidimensional persistence
modules. To the authors' knowledge, such an architecture has not previously been
described, although a special case (where the underlying lattice is a power set)
has been implemented in~\cite{wendler_powerset_2019}. We specialize the convolutions described in
Section~\ref{sec:latticeconv} to the particular case of regular finite
sublattices of $\R^2$. These may be represented (up to isomorphism) as $L = [m]
\times [n]$, where $[n]$ is the ordered set $\{0,1,\dots,n\}$. The meet and join operations are easily computed
elementwise:
\[(r,t) \wedge (r',t') = (\min(r,r'),\min(t,t'));\quad (r,t)\vee (r',t') =
  (\max(r,r'),\max(t,t')).\]

A lattice convolution layer takes as input an $N_{\text{in}}$-dimensional signal
$f : [m] \times [n] \to \R^{N_{\text{in}}}$ and outputs an
$N_{\text{out}}$-dimensional signal $ [m] \times [n] \to
\R^{N_{\text{out}}}$. The layer's parameters are given by a function $g: [m]
\times [n] \to \R^{N_{\text{out}}\times N_{\text{in}}}$. If we label the
entries of $f(x,y)$ by $f_i$ and the entries of $g(x,y)$ by $g^i_j$, the layer then acts by
\[\text{MeetConv}(f)(x,y)^j = \sum_{i} (f_i \ast_{\wedge} g^i_j)(x,y) = \sum_i
  \sum_{(a,b) \in [m]\times [n]} f_i(x \wedge a, y \wedge b)g^i_j(a,b)\]
in the case of convolution with respect to the meet operation, and
\[\text{JoinConv}(f)(x,y)^j = \sum_{i} (f_i \ast_{\vee} g^i_j)(x,y) = \sum_i
  \sum_{(a,b) \in [m]\times [n]} f_i(x \vee a, y \vee b)g^i_j(a,b)\]
in the case of convolution with respect to the join operation.

\paragraph{Remark}
Traditional convolutional neural networks are useful in part because the
convolution kernels (here the functions $g$) can have very small support,
reducing the number of parameters that must be learned. In the standard
convolutional setting, these kernels are implicitly supported in a neighborhood
of the origin, but the location of the kernel is not usually explicitly
specified. In the lattice setting, we do need to specify where the kernel
resides. In the abelian group case, the kernel is supported near the identity,
and similarly, when we treat our domain as a lattice, the kernel should be
supported near the neutral element of the operation. That is, for a meet
convolution, $g$ should be supported at the maximum $(m,n)$, and for a join
convolution, $g$ should be supported at the minimum $(0,0)$. This ensures that
the convolution operators are capable of preserving information at every point
of the space. For instance, if $g(0,0) \neq 0$, then $(f \ast_{\vee} g)(x,y)$ is
a sum of terms including $f(x,y)g(0,0)$, so all information from one layer can
be passed to the next layer. A bit more is necessary to avoid degenerate
convolutions. A kernel supported in a geometrically small neighborhood of the
neutral element results in trivial receptive fields for most neurons. For
example, if $(x,y)$ is greater than the maximal element of the support of $g$,
$(f \ast_{\vee} g)(x,y)$ is simply a scalar multiple of $f(x,y)$. As a result,
we hypothesize that the appropriate support sets for lattice-convolutional kernels
are evenly-spaced sublattices including both $(0,0)$ and $(m,n)$.

\section{Experiments}
We use a portion of the Princeton ModelNet
dataset~\cite{zhirong_wu_3d_2015} as a source of finite
metric spaces. This dataset consists of hundreds of 3-dimensional CAD models
representing objects from 10 classes. We sample
points from the 3D models to produce finite metric spaces embedded in $\R^3$.
We then compute the corresponding multidimensional persistence modules, from
which we produce features used as an input to a convolutional neural net
classifier.

The pipline thus begins with a 3D polyhedral model, of which $3000$ vertices are
sampled to produce a point cloud in $\R^3$. This point cloud then produces a
bifiltered simplicial complex, whose degree-$0$ persistent homology we calculate
using RIVET~\cite{lesnick_interactive_2015}, sampled at a discrete grid of $40
\times 40$ points, producing lattice-indexed signals given by the Hilbert
function and the multi-graded Betti numbers $\xi_0, \xi_1, \xi_2$; four features
in total. These are then passed to the classifier, which produces a class
prediction, in this case one of 10 possible household objects.

As the filter function on these data sets, we use the \introduce{codensity} function \[\rho_{\text{codense}}(x;k) = \left( \frac{1}{k} \sum_{y \in N_k(x)} d(x,y) \right)^{-1},\]
where $N_k(x)$ is the set of the $k$ nearest neighbors to $x$; we select
$k=100$.
\paragraph{Remark}
The name codensity is appropriate because the points in the densest
regions of $\mathcal M$ appear earlier in the filtration. A folk theorem is that
the 2-parameter persistent homology of a Rips/codensity bifiltration is stable (w.r.t.~interleaving distance)
under non-Hausdorff noise (e.g. adding or removing a small number of point samples).

\begin{figure}[h!]
  \caption{The architecture of our lattice-based convolutional neural network.}\label{fig:architecture}
    \begin{center}
    \includegraphics[width=0.75\textwidth]{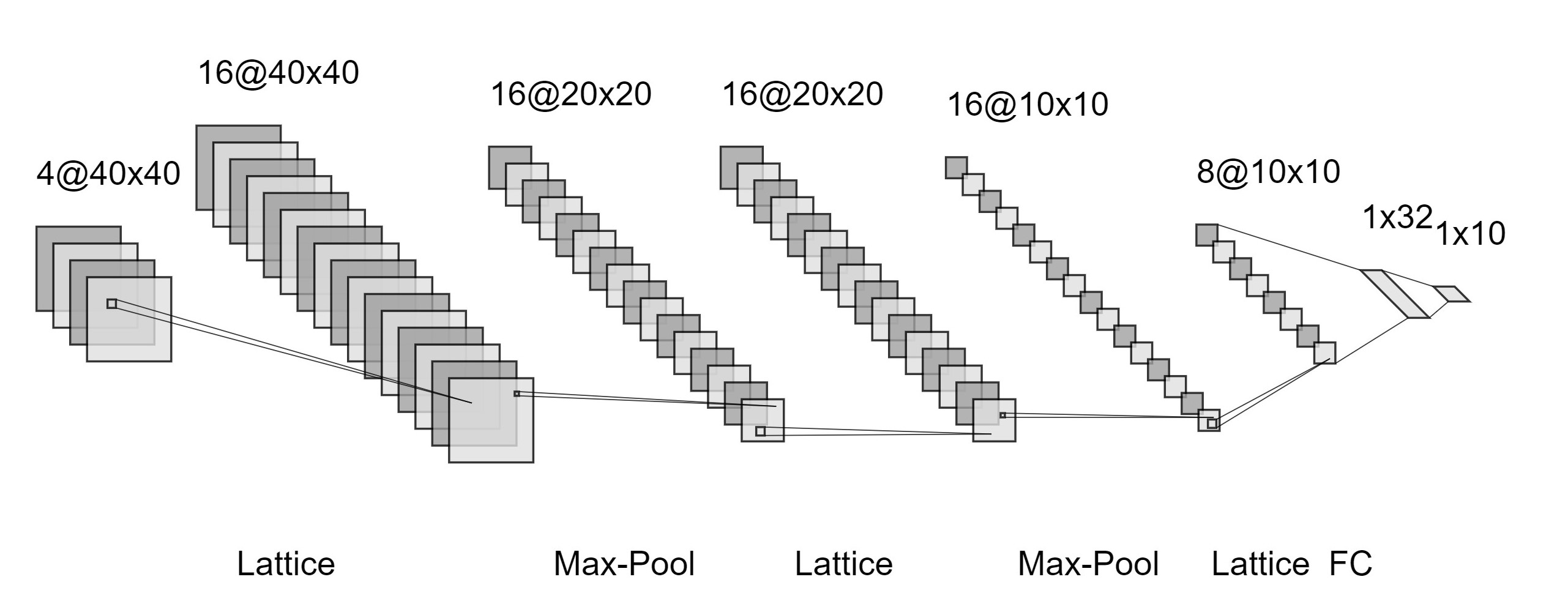}
    \end{center}
\end{figure}

\begin{figure}[h!]
  \caption{Comparison learning curves for the lattice neural network and
    standard convolutional neural network.}\label{fig:comparison}
    \begin{center}
      \includegraphics[width=\textwidth]{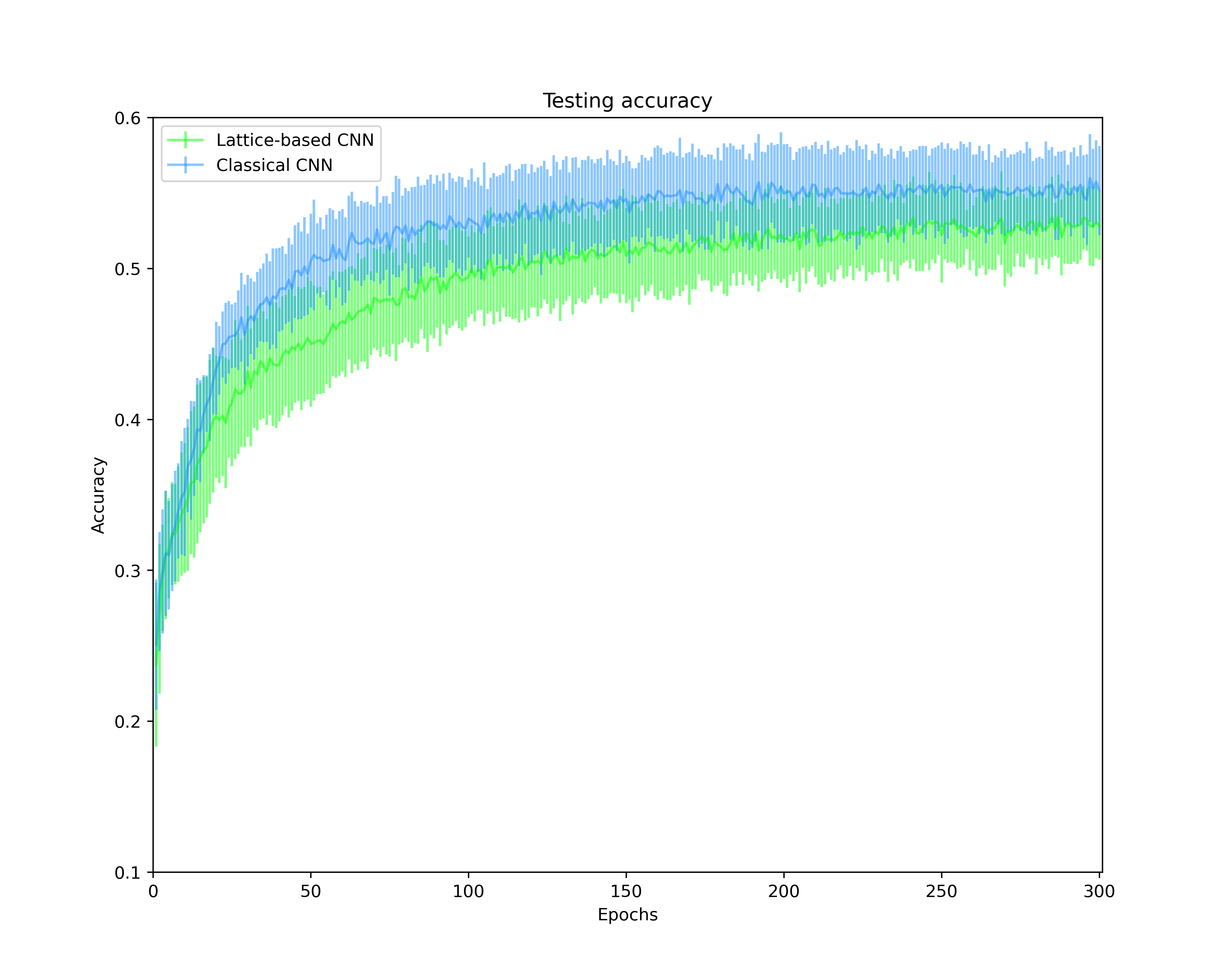}
    \end{center}
\end{figure}

We compare the performance of two convolutional networks on this classification
task. One uses the lattice-convolution based layers described in
Section~\ref{sec:latticeCNN}, and the other uses standard convolutional layers.
Each has three convolutional layers followed by two fully connected layers.
The lattice-based convolution layers are of the form $\alpha \cdot \text{MeetConv}(x) +
(1-\alpha)\cdot\text{JoinConv}(x)$ for a hyperparameter $\alpha \in [0,1]$. We
set $\alpha = \frac{1}{2}$. All
convolution kernels have dimension $4 \times 4$, hidden convolution layers have
16 features, and the final convolution layer has 8 features. The first two convolution
layers are followed by max-pooling layers with a $2 \times 2$ kernel. For the
lattice convolutional layers, the support of the kernel lay in an evenly spaced
$4 \times 4$ grid of points in $[m]\times [n]$, while the standard convolutional
layers had a traditional contiguously supported kernel.
The inner fully connected layer has 32 features. We use a cross entropy loss
function with a softmax in the final layer. The lattice-based convolutional architecture
is summarized in Figure~\ref{fig:architecture}.
The networks are trained using the
Adam gradient algorithm with learning rate $2 \times 10^{-4}$ for 300 epochs. We reserve
10\% of the data for testing. Results are shown in Figure~\ref{fig:comparison}.


The lattice convolutional network slightly underperforms the standard
convolutional classifier. This is somewhat disappointing, and perhaps indicates
that the relevant information encoded in degree-0 persistent homology is best
captured with a geometric notion of locality in $\R^2$.
This work has by no means exhausted the possibilities for lattice convolutions in
multidimensional persistence, but does suggest circumspection in evaluating
their future use.

\section{Discussion}

Our proposed featurization for persistence modules is rather naive, but well
adapted to the use of lattice convolutions as a data processing method. The
lattice convolutional neural network shows promise as a method for classifying
features arising from a multiparameter persistence module. The algebraic
perspective on partially ordered sets, exemplified by lattices, may also offer
approaches to featurizing more complex invariants of persistence modules. In
particular, the incidence algebra may offer a natural way to represent the rank
invariant~\cite{carlsson_theory_2009} in a way amenable to convolution-like
operations.
We hope with these brief experiments to inspire further work on featurizing
multidimensional persistence for use in machine learning algorithms.

\section*{Acknowledgments}
The authors would like to thank Chris Wendler for replicating our original
experiment and pointing out a coding error in our \texttt{pytorch} implementation. We would like to thank the reviewers for providing helpful feedback.
The author [HR] is partially supported by Office of Naval Research (Grant No.
N00014-1442-16-1-2010), and the author [JH] is supported by the National Science
Foundation (NSF-DMS \#1547357).


\end{document}